\documentclass{amsart}
\usepackage{amssymb,amscd,amsthm,verbatim}

\newtheorem{Proposition}{Proposition}
\newtheorem{Theorem}[Proposition]{Theorem}
\newtheorem{Corollary}[Proposition]{Corollary}
\newtheorem{Lemma}[Proposition]{Lemma}

\theoremstyle{definition}

\newtheorem{ex}[Proposition]{Example}

\def\nnull{\underline{\bf 0}}
\def\RR{{\mathbb R}}

\def\rk{{\rm rk}}
\def\ux{{\underline{{\bf x}}}}

\def\Rees{{\mathcal R}}

\newcommand{\NN}{\mathbb N}

\newcommand{\Tor}{\mathrm{Tor}}
\newcommand{\mm}{\mathfrak{m}}

\newcommand{\F}{\mathcal{F}}
\newcommand{\Ga}{\Gamma}
\newcommand{\De}{\Delta}

\numberwithin{equation}{section}

\begin{document}
\title[Segre and Rees products]
{Segre and Rees products of posets, \\
with ring-theoretic applications}

\author[Anders Bj\"orner]{Anders Bj\"orner}
\email{bjorner@math.kth.se}

\address{Department of Mathematics\\
         Royal Institute of Technology\\
         S-10044 Stockholm, Sweden}

\author[Volkmar Welker]{Volkmar Welker}
\email{welker@math.uni-marburg.de}

\address{Fachbereich Mathematik und Informatik \\
         Philipps-Universit\"at Marburg \\
         D-35032 Marburg, Germany}


\thanks{First author was supported by G\"oran Gustafsson Foundation for Research in Natural Sciences and Medicine. Second author was supported by Deutsche Forschungsgemeinschaft (DFG). Both authors' research was supported by EU 
      Research Training Network ``Algebraic Combinatorics in Europe'', grant
      HPRN-CT-2001-00272.}

\begin{abstract}
We introduce (weighted) Segre and Rees products 
for posets and show that these constructions 
preserve the Cohen-Macaulay property over a field $k$ and homotopically.
As an application we show that the weighted Segre product of two affine semigroup rings 
that are Koszul is again Koszul. This result generalizes previous results by  Crona on
weighted Segre products of polynomial rings. 

We also give a new proof of the fact that
the Rees ring of a Koszul affine semigroup ring is again Koszul. 

The paper ends with a list of some open problems in the area.
\end{abstract}

\maketitle

\section{Introduction} \label{intro}

We describe constructions of finite partially ordered sets (posets for short) 
that generalize situations arising in commutative algebra to a
combinatorial setting. For all constructions the principal question 
asked is: ``{\it Does this construction preserve the Cohen-Macaulay property 
?}'' The posets that are relevant for the commutative algebra situation
are those that occur as intervals in affine semigroup posets. 
A result of Peeva, Reiner and Sturmfels \cite{PRS} shows that
an answer to our principal question for the class of 
intervals in affine semigroup posets will give a corresponding 
answer to the question whether certain ring 
theoretic constructions preserve the Koszul property.

The ring-theoretic constructions that motivate this study are
weighted Segre products (see \cite{C})
and Rees algebras.
We define poset-theoretic analogues of these constructions and prove
that the Cohen-Macaulay property is preserved. As corollaries we obtain that weighted 
Segre products of affine semigroup rings preserve the Koszul property. 

\medskip

The following theorem and its corollaries are our main results. Further definitions
and background is given in Section \ref{tools}.
\medskip

\noindent {\sf Segre products of posets:} 
Let $f : P \rightarrow S$ and $g : Q \rightarrow S$ be poset maps.
Let $P \circ_{f,g} Q$ be the induced subposet of the product poset 
$P \times Q$ consisting of the pairs $(p,q) \in P \times Q$ such that 
$f(p) = g(q)$. Recall that the product
poset $P \times Q$ is ordered by $(p,q) \leq (p',q')$ if $p \leq p'$ and
$q \leq q'$. In the language of category theory the poset $P \circ_{f,g} Q$ 
is the pullback of $f : P \rightarrow S \leftarrow Q : g$. 

We will be concerned only with the case when $S = \NN
:= \{0,1,2,\ldots\}$
is the set of natural numbers equipped with its natural order.
For a pure poset $P$ the rank
function serves as an example of a poset map from $P$ to $\NN$.
In case  $P$ is a pure poset with rank function 
$f = \rk$ we write $P \circ_g Q$
for $P \circ_{\rk,g} Q$ and call $P \circ_g Q$ the {\em Segre product}
of $P$ and $Q$ with respect to $g$ (or, the {\em $g$-weighted Segre product }
of $P$ and $Q$).

\begin{Theorem} \label{poset-segre}
Let $P$ and $Q$ be pure posets. Let $\rk : P \rightarrow \NN$ be the rank
function of $P$ and $g : Q \rightarrow \NN$ be a strict poset map such that
$g(Q) \subseteq \rk (P)$. If $P$ and $Q$ are Cohen-Macaulay over the field $k$,
then the Segre product $P \circ_{g} Q$ is Cohen-Macaulay over $k$.
If $P$ and $Q$ are homotopically Cohen-Macaulay, then so is $P \circ_{g} Q$.
\end{Theorem}

Theorem \ref{poset-segre} is a special case of a more general result for simplicial complexes,
see Theorem \ref{complex-segre} below.
\medskip

\noindent {\sf Rees products of posets:}
Let $P$ and $Q$ be pure posets with rank functions $\rk$. Let $P*Q$ be the 
poset on the ground set $\{ (p,q)\in P\times Q~|~\rk(p) \geq \rk(q) \}$ with 
order relation 
$$(p,q) \leq (p',q') \stackrel{\text{def}}{\Longleftrightarrow}  
p \leq p' , q \leq q' \mbox{~and~} \rk (p') - \rk (p) \geq \rk(q') - \rk(q).$$
We call $P*Q$ the {\em Rees product} of $P$ and $Q$. Note that it is not in
general an induced subposet of the product $P \times Q$. However,
as will be shown, it is nevertheless a special case of the Segre product. 
Thus, using
Theorem \ref{poset-segre} we prove

\begin{Corollary} \label{poset-rees}
Let $P$ and $Q$ be pure posets. 
If $P$ and $Q$ are Cohen-Macaulay over the field $k$
and $Q$ is acyclic over $k$,
then the Rees product $P * Q$ is Cohen-Macaulay over $k$.
If $P$ and $Q$ are homotopically Cohen-Macaulay and $Q$ is
contractible, then $P * Q$ is homotopically Cohen-Macaulay.
\end{Corollary}

\medskip

\noindent {\sf Affine semigroup rings:} 
Theorem \ref{poset-segre} has the following ring-theoretic consequence,
explained and further discussed in Section \ref{affinerings}.

\begin{Corollary} \label{ring-segre}
Let $\Lambda \subseteq \NN^d$ and $\Gamma \subseteq \NN^e$
be two homogeneous affine semigroups. Assume that the semigroup rings 
$k[\Lambda]$ and
$k[\Gamma]$ are Koszul. Let $g$ be a grading of  $k[\Gamma]$. 
Then the weighted Segre product $k[\Lambda] 
\circ_g k[\Gamma]$ with respect to $g$ is Koszul.
\end{Corollary}

Similarly, Theorem \ref{poset-rees} has the following consequence,
which can also be deduced by ring-theoretic arguments from a result
of Backelin and Fr\"oberg \cite[Proposition 3]{BF}.

\begin{Corollary} \label{ring-rees}
Let $k[\Lambda]$ and $k[\Gamma]$ be Koszul affine semigroup rings. Let
$k[\Lambda]_i$ and $k[\Gamma]_i$ denote their $i$-th graded components and
set $\mm_\Lambda = \bigoplus_{i \geq 1} k[\Lambda]_i$.
Then the $k$-algebra
$$k[\Lambda] * k[\Gamma] = \bigoplus_{i \geq 0} \mm_\Lambda^i 
\otimes_k k[\Gamma]_i$$ is Koszul.
\end{Corollary}

Note that for $k[\Gamma] = k[t]$ the Rees product 
$k[\Lambda] * k[\Gamma]$ is the Rees ring $\Rees[k[\Lambda],\mm_\Lambda]$
of $k[\Lambda]$ with respect to its maximal ideal $\mm_\Lambda$.

\section{Tools from Topological Combinatorics}\label{tools}

We begin with a review of some basic definitions.

A {\em chain} $C$ in a poset $P$ is a linearly ordered subset,
its {\em length} $\ell (C)$ is one less than its number of elements.
A  poset $P$ is {\em pure} if all maximal chains have the
same length. 
For each element $p \in P$ of a pure poset $P$ the length of a maximal
chain in $P_{\le p} := \{q\in P\mid q\le p\}$ is
called the {\em rank} $\rk(p)$ of $p$ in $P$.

A poset $P$ is called {\em bounded} if there is a unique minimal element 
$\hat{0}$ and a unique maximal element $\hat{1}$ in $P$. 
For two elements $x \leq y$ in $P$ we write $[x,y]$ for the {\em closed interval}
$\{ z~|~x \leq z \leq y \}$ in $P$, and similarly
$(x,y)$ for the {\em open interval}
$\{ z~|~x < z < y \}$.
Clearly, $[x,y]$ is a bounded poset.
A poset $P$ is called {\em graded} if it is both
bounded and pure. Let $\widehat{P} :=P\cup\{\hat{0},\hat{1}\}$
denote $P$ augmented by new bottom and top elements
$\hat{0}$ and $\hat{1}$. Thus, $\widehat{P}$ is graded
iff $P$ is pure. 

A map $f:P\rightarrow Q$ is a {\em poset map} 
[resp.  a {\em strict poset map}]
if $x<y$ implies $f(x)\le f(y)$ [resp. $f(x)<f(y)$] for all $x,y\in P$.

We write $\Delta(P)$ for the simplicial
complex of all chains of $P$. 
By $\widetilde{H}_i(P;k)$ we denote the $i$-th reduced simplicial 
homology group 
of $\Delta(P)$ with coefficients in $k$. 
Also, if convenient we identify $P$ with the geometric realization of
$\Delta(P)$. 

A poset $P$ is called 
{\it Cohen-Macaulay} over the field $k$ if for all $x < y$ in $\widehat{P}$
the reduced simplicial homology $\widetilde{H}_i((x,y)\,;k)$ vanishes for 
$i \neq \rk(y) - \rk(x) - 2$. A poset $P$ is called {\em homotopically} 
Cohen-Macaulay if for all $x < y$ in $\widehat{P}$ the interval $(x,y)$ is homotopy equivalent
to a wedge of spheres of dimension $\rk(y) - \rk (x) - 2$.
Cohen-Macaulay posets are pure.

The Cohen-Macaulay-over-a-field-$k$ property for posets is a special case of a
property defined for all finite simplicial complexes. This general notion
of Cohen-Macaulayness is in turn equivalent to a 
particular instance of the ring-theoretic Cohen-Macaulay concept.
For this connection with Commutative Algebra, see Stanley \cite{St}.
\medskip

The construction of Segre products of posets generalizes a well
known concept; namely, Segre products of posets subsume 
rank selection of posets. 
For a pure poset $P$ with rank function $\rk$,
let $S \subseteq \rk(P)$ be a set of ranks of the poset $P$.
The {\em rank-selected subposet} of $P$ determined by $S$ is the 
induced subposet 
$P_S$ of all elements $x$ of $P$ such that $\rk(x) \in S$.

In the late 1970's, rank selection was shown to preserve
Cohen-Macaulayness  over $k$ by Baclawski, Munkres, Stanley and
Walker. See e.g.  \cite{B}. Complete references are given in \cite[p. 1858]{Bj},
where also a proof of the homotopy version is sketched.

\begin{Proposition}\cite[Thm. 6.4]{B} \cite[Thm. 11.13]{Bj} \label{rank}
\label{rankselection}
Let $P$ be a Cohen-Macaulay poset over the field $k$ and let
$S \subseteq \rk(P)$.
Then the rank-selected subposet $P_S$ is Cohen-Macaulay over $k$.
If $P$ is homotopically Cohen-Macaulay, then so is $P_S$.
\end{Proposition}

In order to realize rank selection as a Segre
product, let $Q$ be the chain on $|S|$ elements $\{ 1,
\ldots,r\}$. Let $S = \{ s_1 < \cdots < s_r \}$ and let
$g$ be the map that sends $i \in Q$ to $s_i$. Then 
it is easily seen that $P \circ_g Q \cong P_S$. Unfortunately,
Theorem \ref{poset-segre} does not give a new proof of 
Proposition \ref{rankselection}, but rather 
uses this fact as an essential point in the argumentation.

\bigskip

As a second tool we need another result, which is due to Baclawski \cite{B}
for Cohen-Macaulayness over $k$ and to Quillen \cite{Q} for
homotopical Cohen-Macaulayness. Both versions are proved in slightly
greater generality  in \cite{BWW}.

\begin{Proposition} \cite[Thm. 5.2]{B},  \cite[Cor. 9.7]{Q}
 \label{fibration} 
Let $P$ and $Q$ be pure posets and $f : P \rightarrow Q$ a rank-preserving and 
surjective poset map. 
Assume that for all $q \in Q$ the fiber $f^{-1}(Q_{\leq q})$ is
Cohen-Macaulay over $k$. If
$Q$ is Cohen-Macaulay over $k$, then so is also $P$.
The same is true with ``Cohen-Macaulay over $k$'' everywhere replaced by 
``homotopically Cohen-Macaulay''.
\end{Proposition}

We also need  the following result
on barycentric subdivisions, which can be obtained from the
fact that Cohen-Macaulayness is invariant under homeomorphisms.

First recall a few definitions. 
If $\Gamma$ is a simplicial complex then we can consider $\Gamma$ as a
poset, namely as the
partially ordered set of its faces ordered by inclusion.
Denote by $\F (\Ga )$ this {\em face poset}.
The simplicial complex 
$\Delta(\F(\Gamma))$ is called the
{\it barycentric subdivision} of $\Gamma$ and is well know to be
homeomorphic to $\Gamma$.  

\begin{Proposition} \label{barycenter} 
A poset $P$ is Cohen-Macaulay over $k$ (resp. homotopically
Cohen-Macaulay) if and only if the poset $\F(\De(P))$ has the
same property.
\end{Proposition}

\section{Proofs and comments}\label{proofs}

In this section we prove the main poset theoretic theorems
and discuss some related questions.

\medskip

The $g$-weighted Segre product  $P \circ_{g} Q$ of two pure posets $P$ and $Q$
was defined in Section \ref{intro}. Note that $P \circ_{g} Q$ is also pure, and that
$\rk_{P \circ_{g} Q} (p,q) =\rk_{Q}(q)$, for all $(p,q)\in P \circ_{g} Q$. 
In particular, $\rk(P \circ_{g} Q)=\rk(Q)$.
\medskip

\noindent {\sf Proof of Theorem  \ref{poset-segre}.}
Let $$f: \F(\De (P \circ_{g} Q ))\rightarrow \F(\De (Q ))$$ be the poset map that sends
each chain $(p_0,q_0) < \cdots < (p_{\ell} ,q_{\ell})$ to its projection 
$q_0 < \cdots < q_{\ell}$. This map is surjective and rank-preserving.
For an element $c = (q_0 < \cdots < q_{\ell})$ in $\F(\Delta(Q))$,
the fiber $f^{-1}(\F({\Delta(Q)})_{\leq c})$ consists of all subchains of chains
$(p_0,q_0) < \cdots < (p_{\ell},q_{\ell})$ for which $\rk(p_i) = g(q_i)$ for all $i$. 
Setting $S = \{ g(q_0), \ldots, g(q_{\ell})\}$, then clearly
$f^{-1}(\F({\Delta(Q)})_{\leq c})$
is isomorphic to $\F(\Delta(P_S))$. 

By Proposition \ref{rankselection} we know that $P_S$ is 
Cohen-Macaulay over $k$ (resp. homotopically Cohen-Macaulay), since $P$ is. Also,
Proposition \ref{barycenter} shows that since $P_S$ is  Cohen-Macaulay
over $k$ (resp. homotopically Cohen-Macaulay) then so is also 
$\F(\Delta(P_S))$. Hence, we get from Proposition \ref{fibration} 
that $\F(\De (P \circ_{g} Q ))$
is Cohen-Macaulay. The assertion now follows 
via Proposition \ref{barycenter}. \hfill
$\Box$
\medskip

We don't see any reasonable way to go beyond  Theorem \ref{poset-segre}
in its poset version. Consider these obstacles:
\begin{itemize}
\item If $g(Q) \subseteq \rk(P)$ is not required, then if $P$ is a chain we can
realize arbitrary lower order ideals in $Q$ as Segre products $P \circ_g Q$.
\item If $g$ is not {\it strict},  then a counterexample to the conclusion of
the theorem can be constructed as follows. Let $P$ 
be a two element antichain, let $Q = \{x < y \}$ be a
two element chain, and let $g(x) = g(y) = 0$.
Then $P \circ_g Q$ is the disjoint union of two chains of
length $1$, and hence the poset is not Cohen-Macaulay.
\end{itemize}
However, there is a rather 
straight-forward generalization of the Segre product  to simplicial complexes.

Let $\Ga_1$ and $\Ga_2$ be simplicial complexes on vertex sets
$V_1$ resp. $V_2$, with $\dim \Ga_2 \le \dim \Ga_1 =d-1$.
Assume that there are maps $g_i : V_i \rightarrow \{1,\ldots ,d\}$ such that:
\begin{enumerate}
\item[(i)] $g_1$ restricts to a bijection on each maximal face of $\Ga_1$, 
\item[(ii)] $g_2$ is injective on each maximal face of $\Ga_2$.
\end{enumerate}
Define a simplicial complex $\Ga_1 \circ_{g_1, g_2} \Ga_2$ on 
the vertex set $V_1 \times V_2$ as having faces
$$\{(x_1, y_1), \ldots, (x_k, y_k) \}$$
for all $\{x_1, \ldots, x_k \}\in\Ga_1$ and $\{y_1, \ldots, y_k \}\in\Ga_2$
such that $g_1(x_i)=g_2(y_i)$ for all $i$.
\begin{Theorem}\label{complex-segre}
If \,$\Ga_1$ and $\Ga_2$ are Cohen-Macaulay over $k$ (resp. homotopically
Cohen-Macaulay), then so is $\Ga_1 \circ_{g_1, g_2} \Ga_2$.
\end{Theorem}
\begin{proof}
Essentially the same proof as for Theorem \ref{poset-segre} goes through.
Instead of ``rank-selected subposets'' one has here to use 
``type-selected subcomplexes'', for which the preservation of
Cohen-Macaulayness is also known, see \cite[p. 1858]{Bj}.
\end{proof}

In our opinion, even the following specialization of Theorem \ref{poset-segre},
to what might be called ``unmixed Segre products'',
is somewhat unexpected from the combinatorial point of view.

\begin{Corollary} \label{usual-poset-segre} Let $P$ and $Q$ be
pure posets and let $\rk$ denote the rank function for either poset. 
If $P$ and $Q$ are Cohen-Macaulay over $k$ then the poset 
$P \circ_\rk Q = \{ (p,q)~|~\rk (p) = \rk (q) \}$ is Cohen-Macaulay
over $k$. 
If $P$ and $Q$ are homotopically Cohen-Macaulay, then so is $P \circ_{\rk} Q$.  
\end{Corollary}

\begin{ex}
Let $M_n$ denote the poset of all minors (square submatrices)
of an $n\times n$ matrix.  As a special case of Corollary \ref{usual-poset-segre} 
one sees that this {\em poset of minors} is Cohen-Macaulay. Namely,
if $B_n$ denotes the Boolean lattice of all subsets of $[n] := \{1,\ldots, n\}$,
then $M_n$ is clearly isomorphic to the {\em Segre square}
$B_n \circ_\rk B_n =\{(A,B)\mid A,B\subseteq [n], |A|=|B|\} 
\subseteq B_n\times B_n$. 
Such Segre powers (of infinite posets) 
previously appeared in the work of Stanley, see \cite[Example 1.2]{St2}.

The number of $(n-2)$-spheres in the wedge giving the
homotopy type of $\De(M_n \setminus \{\hat{0}, \hat{1}\})$,
or equivalently $(-1)^{n} \mu( \hat{0}, \hat{1})$ 
where $\mu( \hat{0}, \hat{1})$
is the value of the M\"obius function over $M_n$, is equal to
the number of
pairs of permutations of $[n]$ having no common ascent. This
set of permutation-pairs is well studied (see \cite{CSV}). In \cite{CSV}
one can find a recurrence relation for these numbers which is
exactly the defining relation for the M\"obius number of the 
Segre square of $B_n$.

A second way to obtain this enumerative result is via the theory
of lexicographic shellability \cite{BjWa2}. A natural labeling rule for
$B_n \circ_\rk B_n$ is to give a covering
$(A_1 , B_1) \subset (A_2 , B_2)$
the label $(a,b)$, where
$a$ and $b$ are the unique elements of $A_2 - A_1$
and $B_2 - B_1$, respectively. This is clearly an EL-labeling,
and the falling chains are labeled by 
pairs of permutations with no common ascent.

A third approach is via the rank-selected
$\alpha$- and $\beta$-invariants 
$\alpha_J$ and $ \beta_J$ of $B_n$, as defined by Stanley \cite[p.131]{St1}.
One gets that
$$(-1)^{n} \mu( \hat{0}, \hat{1})=\sum_{J\subseteq [n-1]} \alpha_J \beta_J.$$
This expression for $\mu( \hat{0}, \hat{1})$ of $B_n \circ_\rk B_n$
follows from
\cite[Theorem 5.1 (iii)]{BWW}, and is more generally true for 
Segre squares of all
Gorenstein*  (i.e., Cohen-Macaulay and Eulerian) posets.
Since the Boolean lattice is lexicographically shellable there is a simple 
interpretation of $\alpha_J$ and $\beta_J$. For lexicographically 
shellable posets $\beta_J$ counts the number of maximal chains whose
descent set is equal to $J$ and $\alpha_J$ counts the number of 
maximal chains whose descent set is contained in $J$.
If one uses the labeling $\lambda$
of cover relations in $B_n$ where $\lambda(A \subset B)$ is the
unique element of $B - A$, then 
maximal chains correspond to permutations in $S_n$. Thus 
$\displaystyle{\sum_{J \subseteq [n-1]}} \alpha_J \beta_J$
counts pairs of permutations  $(\sigma, \tau)$ such that the descent 
set of the first is contained in the descent set of the second. 
Equivalently, it counts pairs $(\sigma, \tau)$ of permutations
such that at a place where $\sigma$ has a
descent the permutation $\tau$ has an ascent. Now, if we reverse 
the permutations 
(when written as words) this set bijects
to pairs of permutations with no common ascent.
\end{ex}
\bigskip

We now turn to the
Rees product  $P * Q$ of two pure posets $P$ and $Q$,
defined in Section \ref{intro}. Note that $P * Q$ is also pure, and that
$\rk_{P * Q} (p,q) =\rk_{P}(p)$, for all $(p,q)\in P * Q$. 
In particular, $\rk(P *Q)=\rk(P)$.

\begin{Lemma}
Let $P$ and $Q$ be pure posets. Furthermore, let $\widetilde{Q}:= (Q \times C_n)_{[0,n]}$,
where $n = \rk (P)$, $C_n$ is 
a chain of $n+1$ elements,  and the subscript denotes rank-selection
to the elements of rank at most $n$ in the direct product.
Then the Rees product $P * Q$ is
isomorphic to the (unweighted) Segre product $P \circ \widetilde{Q}$.
\end{Lemma}
\begin{proof}
The elements of $P \circ \widetilde{Q}$ are of the form
$(p,q,i)$, where $\rk(p) = \rk (q) + i$, $0\le i\le n$. In particular, $\rk(q) \leq \rk(p)$.
Now we have $(p,q,i) \leq (p',q',i')$ if and only if $p \leq p'$, $q \leq
q'$ and $i \leq i'$. Thus by $\rk(p) = \rk (q) + i$ and 
$\rk(p') = \rk (q') + i'$ we infer that $\rk(p) - \rk(q) \leq \rk(p')
- \rk(q')$. Thus the projection map onto the first two
coordinates is an isomorphism from $P \circ \widetilde{Q}$ to
the Rees product $P * Q$.
\end{proof}   

\noindent {\sf Proof of Corollary \ref{poset-rees}.}
Let $C_n$ denote a chain of $n+1$ elements,
where $n = \rk (P)$.
By results
of Baclawski \cite{B} and Walker \cite{Wa}, a direct product of two posets which
are Cohen-Macaulay over $k$ (resp. homotopically Cohen-Macaulay) is again
Cohen-Macaulay over $k$ (resp. homotopically Cohen-Macaulay) 
if both posets are acyclic over $k$ (resp. contractible).
Thus $Q \times C_n$ is Cohen-Macaulay if $Q$ is Cohen-Macaulay and 
acyclic over $k$ (resp. contractible). By Proposition \ref{rank} then
also $\widetilde{Q}$ is Cohen-Macaulay, and finally it follows from
Theorem \ref{poset-segre}  that $P \circ \widetilde{Q} \cong P*Q$ is
Cohen-Macaulay. 
\hfill $\Box$  

\begin{ex}\label{deranged}
Let $B_n \setminus \{ \emptyset \}$ be the Boolean lattice of all subsets of 
$[n]$ with the empty set removed. Let $C_n$ be a chain of $n$ elements.
Both $B_n \setminus \{ \emptyset \}$ 
and $C_n$ are homotopically Cohen-Macaulay
and contractible.
By Corollary \ref{poset-rees} we therefore know that
$R_n= (B_n \setminus \{ \emptyset \}) * C_n$ is 
homotopically Cohen-Macaulay.

Attempts to compute the exact homotopy type of the poset $R_n$ have 
led to a problem that we state at the end of Section 5.

\end{ex}

\newcommand{\ME}{{\mathcal{ME}(Q)}}
\newcommand{\M}{{\mathcal{M}(Q)}}
\newcommand{\la}{\lambda}
\newcommand{\cov}{\triangleleft}

\section{Affine semigroup rings}\label{affinerings}

Our initial motivation for this work comes 
from the study of the Koszul property for 
affine semigroup rings in commutative algebra. In this section we explain this
motivation and the ring-theoretic consequences of our main results.

\smallskip

Let $\Lambda \subseteq \NN^d$ be an affine semigroup 
(i.e., a finitely generated additive sub-semigroup containing $\nnull$). For 
$\lambda = (\lambda_1,\ldots,\lambda_d)\in \Lambda $
we set $\ux^{\lambda} = x_1^{\lambda_1} \cdots x_d^{\lambda_d}$.
For a field $k$
the semigroup-ring $k[\Lambda] \subseteq k[x_1,\ldots,x_d]$ is the subalgebra 
of the polynomial ring $k[\NN^d] = k[x_1,\ldots,x_d]$
generated by all monomials $\ux^\lambda$ for
$\lambda \in \Lambda$. 

The semigroup $\Lambda$ is equipped with the
structure of a partially ordered set by setting $\lambda \leq \gamma$
if there is a $\rho \in \Lambda$ such that $\lambda + \rho = \gamma$.
Clearly, $\nnull$ is the unique minimal element of the semigroup $\Lambda$
regarded as a poset. In the sequel we will always assume that
the elements of $\Lambda$ span $\RR^d$ as a vector space. Then as posets
all intervals in $\Lambda$ are pure if all elements of a minimal generating 
set lie on an affine
hyperplane; in this situation we also say $\Lambda$ is {\em homogeneous}.
Note, that by the commutativity of $\Lambda$ it follows that 
every lower interval $[\nnull,\lambda]$ is self-dual as a poset.
Also, in the poset $\Lambda$ every interval is isomorphic to a lower
interval: $[\mu,\lambda] \cong [\nnull,\lambda-\mu]$.

We call a $k$-algebra 
$A$ {\it standard graded} if as a $k$-vector space $A \cong 
\bigoplus_{i \in \NN} A_i$, $A_0 = k$, $A_iA_j \subseteq A_{i+j}$ and 
$A$ is as an algebra generated by $A_1$. 
If $\Lambda$ is homogeneous then
$k[\Lambda]$ is a standard graded algebra. 

A standard graded $k$-algebra $A = \bigoplus_{i \in \NN} A_i$ is called 
{\it Koszul} if $k$ has a linear resolution over $A$;
or equivalently, if $\Tor_i^A(k,k)_j =0$ for $i \neq j$
(see \cite{F} for a comprehensive survey on 
Koszul rings). Via the bar resolution
and work of Laudal and Sletsj\o e \cite{LS},
Peeva, Reiner and Sturmfels \cite{PRS} observe the following relation between
the Koszul property and Cohen-Macaulayness for affine semigroup rings:

\begin{Proposition}[\cite{PRS}] \label{CM-Koszul} For an affine semigroup
$\Lambda$ and a field $k$ the following are equivalent:
\begin{itemize}
\item[(i)] the ring $k[\Lambda]$ is Koszul;
\item[(ii)] the interval $(\nnull,\lambda)$ is 
a Cohen-Macaulay poset over $k$, for all $\lambda \in \Lambda$;
\item[(iii)] the interval $(\nnull,\lambda)$ 
is pure and has homology concentrated in dimension $\rk(\lambda) -2$,
for all $\lambda \in \Lambda$.
\end{itemize}
\end{Proposition}

Using this lemma, we now draw the ring-theoretic conclusions of
our work in earlier section. For this we first review
the required ring-theoretic concepts.
\medskip

\noindent {\sf Weighted Segre products:}
Let $A = \bigoplus_{i \geq 0} A_i$ and $B = \bigoplus_{i \geq 0} B_i$ be two
graded $k$-algebras. Also, let $B = \bigoplus_{i \geq 0} B_i'$ be another 
grading of $B$ as a $k$-algebra; i.e., $B'_iB'_j \subseteq B'_{i+j}$. 
Assume that as $k$-vector spaces $B_i = \bigoplus_{j \geq 0} 
B_i \cap B_j'$. The {\it weighted Segre product} 
$A \circ' B$ of $A$ and $B$ with respect to the grading $B = \bigoplus_{i \geq 0} B_i'$
is the $k$-subalgebra of $A \otimes B$ 
generated by the elements $a \otimes b \in A \otimes B$ such that 
$a \in A_i$ and $b \in B_i'$ for $i \in \NN$.
If $A$ and $B$ are standard graded $k$-algebras 
then $A \circ' B$ is generated
as a $k$-algebra by $\bigoplus_{i \geq 0} A_i \otimes (B_i' \cap B_1)$. 
The concept of a weighted Segre product first appeared in work of Crona
\cite{C}, where weighted Segre products of polynomial rings are considered. 
Since our results on weighted Segre products 
apply to affine semigroup rings only, we from now on we confine ourselves
to this setting. \\

Let $k[\Lambda]$ and $k[\Gamma]$ be two affine semigroup rings
for the homogeneous affine semigroups $\Lambda \subseteq \NN^d$ and 
$\Gamma \subseteq \NN^e$. Let $f$ be the standard grading for $\Lambda$. 
Let $g : \Gamma \rightarrow \NN$ be some grading (i.e., semigroup
map with $g(\gamma)>0$ for all $\gamma\neq 0$).\\ 

The {\it weighted Segre product} of the affine
semigroups $\Lambda$, $\Gamma$, with respect to the grading $g$,
 is the affine semigroup 
$\Lambda \circ_g \Gamma \subseteq \NN^{d+e}$ of 
all pairs $(\lambda,\gamma)$ with $f(\lambda) = g(\gamma)$.
One easily sees that the semigroup-ring  $k[\Lambda \circ_g \Gamma]
\subseteq k[\NN^{d+e}]$ is (isomorphic to) the weighted Segre product
$k[\Lambda] \circ_g k[\Gamma]$
(in the sense of the previous paragraph) of the 
affine semigroup rings $k[\Lambda]$ and $k[\Gamma]$.
The semigroup ring $k[\Lambda] \circ_g k[\Gamma]$ is again
homogeneous with grading induced by $(\lambda,\gamma) \mapsto h(\gamma)$, where
$h$ is the standard grading for $\Gamma$.

\medskip

We can now derive
Corollary \ref{ring-segre} from Theorem \ref{poset-segre}.

\noindent {\sf Proof of Corollary  \ref{ring-segre}.}
It is easily seen that 
if $(\lambda,\gamma) \in \Lambda \circ_g \Gamma$ then the lower interval
$[\nnull, (\lambda,\gamma)]$ in $\Lambda \circ_g \Gamma$ is isomorphic
to the $g$-weighted Segre product of posets $[\nnull,\lambda] \circ_g 
[\nnull,\gamma]$. Hence, Theorem \ref{poset-segre} implies
Corollary \ref{ring-segre} via Proposition \ref{CM-Koszul}. 
\hfill$\Box$
\medskip

We describe some special cases. 

\begin{itemize}
\item {\tt Segre product}: 
      If $g$ is the standard grading of $k[\Gamma]$ then 
      $k[\Lambda] \circ_g 
      k[\Gamma]$ is the usual Segre product 
      $k[\Lambda] \circ k[\Gamma]$ of rings. It is known that 
      in general the Segre product of two Koszul rings is again Koszul 
      (Backelin \& Fr\"oberg \cite{BF}).  
\item {\tt Veronese-ring}: 
      If $\Gamma = \NN$ and $g(1) = s$ then $k[\Lambda] \circ_g
      k[\Gamma]$ is the $s$-th Veronese ring of $k[\Lambda]$.
      Again it is known that in general a Veronese ring of a 
      Koszul ring is Koszul (Backelin \& Fr\"oberg \cite{BF}).  
\item {\tt Polynomial rings}: If $\Lambda = \NN^d$ and $\Gamma =
      \NN^e$ then for a grading $g : \Gamma
      \rightarrow \NN$ such that $\Gamma$ is generated in a fixed 
      $g$-degree the ring $k[\Lambda] \circ_g k[\Gamma]$
      is Koszul (Crona \cite{C}).
\end{itemize}
\medskip

\noindent {\sf Rees products:} Let $A$ be a ring and $I$ an ideal in $A$.
Then the {\em Rees ring} $\Rees[A,I]$ is the direct sum
$\bigoplus_{i \geq 0} t^i I^i$, where $t$ is an additional indeterminate and
$I^0 = A$. Here we consider the case when $A = \bigoplus_{i \geq 0} A_i$
is a standard graded $k$-algebra and $I = \mm_A = \bigoplus_{i \geq 1} A_i$.
We also generalize the construction in the following way. Let $B = \bigoplus_{i \geq 0} B_i$
be another standard graded $k$-algebra. Then we define the {\em Rees product} 
$A * B$ as the $k$-algebra $\bigoplus_{i\geq 0} \mm_A^i \otimes_k B_i$. 
If $B = k[t]$ is the polynomial ring in a single variable the Rees product
$A * B$ is the Rees ring $\Rees[A,\mm_A]$.

Essentially the same arguments that show that Rees rings of a Koszul algebra
with respect to the maximal ideal are Koszul also show that the Rees product
$A * B$ preserves Koszulness.

\begin{Proposition} \label{gen-ring-rees} Let $A$ and $B$ be Koszul standard graded $k$-algebras.
Then $A * B$ is Koszul.
\end{Proposition}
\begin{proof}
Consider the Segre product $R = A \circ (B \otimes k[t])$. It is easily seen that the
projection on $A * B$ is a $k$-algebra isomorphism.
Moreover, by \cite{BF} we know that Segre products preserve Koszulness, as do
tensor products. Thus $R$ is a Koszul $k$-algebra.
\end{proof}

We consider the case when $A = k [\Lambda]$ and $B = k[\Gamma]$ are
standard graded affine semigroup rings for semigroups $\Lambda \subseteq \NN^d$ 
and $\Gamma \in \NN^e$.
One checks that $k[\Lambda] * k [\Gamma]$ is the affine semigroup ring
$k[\Lambda * \Gamma]$, where $\Lambda * \Gamma \subseteq \NN^{d+e}$ is the
affine semigroup generated by $(\lambda,\nnull)$ and $(\lambda,\gamma)$ for
elements $\lambda \in \Lambda$ and $\gamma \in \Gamma$ of degree $1$.
Clearly, Proposition \ref{gen-ring-rees} implies that for Koszul $k[\Lambda]$ and
$k[\Gamma]$ the Rees product $k[\Lambda] * k[\Gamma]$ is Koszul as well.
But we want to present an alternative derivation of this fact by using
the poset Rees product in order to give the motivation for our poset theoretic
construction.  

\medskip

\noindent {\sf Proof of Corollary  \ref{ring-rees}.}
Let $\rk_\Lambda$ and $\rk_\Gamma$ be the rank functions of
$\Lambda$ and $\Gamma$, and 
let $(\lambda', \gamma') \in \Lambda * \Gamma$. Since 
$\Lambda * \Gamma$ is generated be elements $(\lambda,\gamma)$ where
$\rk_\Lambda \lambda \geq \rk_\Gamma \gamma$, it follows that $\rk_\Lambda
\lambda'' \geq \rk_\Gamma \gamma''$ for all $(\lambda'', \gamma'')
\leq (\lambda',\gamma')$. Moreover, this also implies that 
$(\lambda'', \gamma'') \leq (\lambda',\gamma')$ if and only if 
$\rk_\Lambda \lambda' -\rk_\Lambda \lambda'' \geq
\rk_\Gamma \gamma' -\rk_\Gamma \gamma''$.
Thus $$[\nnull,(\lambda', \gamma')] \cong [\nnull,\lambda] *
[\nnull,\gamma].$$
Now, by Proposition \ref{CM-Koszul} we get that $[\nnull,\lambda]$
and $[\nnull,\gamma]$ are Cohen-Macaulay over $k$. Having
a least and a maximal elements implies that $[\nnull,\gamma]$ is
contractible. Thus by
Corollary \ref{poset-rees} it follows that $[\nnull,\lambda] *
[\nnull,\gamma]$ is Cohen-Macaulay over $k$. Another application
of Proposition \ref{CM-Koszul} then proves the assertion.
\hfill$\Box$

\medskip

\section{Some open problems}\label{problems}

The connection between topological combinatorics and ring theory 
via semigroup posets offers several interesting open problems.
In closing we list a few.
\smallskip

In the following $\Lambda\subseteq \NN^d$ 
is a homogeneous affine semigroup, ordered in the 
usual way. Let $\rk$ be the rank function of $\Lambda$ as a poset. Then the
dimension $\dim((\nnull,\la))$ of the order complex of the open interval $(\nnull,\la)$ is $\rk (\la)-2$.

Let us call $\Lambda$ Cohen-Macaulay over a
fixed field $k$ (resp. homotopically Cohen-Macaulay), if all intervals 
$(\nnull, \la)$ in $\Lambda$ are Cohen-Macaulay over $k$
(resp. homotopically Cohen-Macaulay). Clearly, Cohen-Macaulayness 
over $k$ depends on the characteristic of the field $k$ only,
and this property is equivalent to Koszulness of the ring
$k[\Lambda]$ (by Proposition \ref{CM-Koszul}).

\begin{enumerate}

\item[(1)] 
Ring-theoretic work of Avramov and Peeva \cite{AP} 
implies the following fact:
Suppose that there exists a $\la\in\Lambda$ and an $i>0$ such that
$$\tilde{H}_{\rk(\la)-2-i}((\nnull,\la)\,;k)\neq 0.$$
Then for all $j>0$ there exists some $\la'\in\Lambda$ and some
$j' \geq j$ such that
$$\tilde{H}_{\rk(\la')-2-j'}((\nnull,\la')\,;k)\neq 0.$$ 

\noindent {\sf Question:} Does this have a combinatorial explanation?

Moreover, given $\lambda$, $i$ and $j$ it would be interesting to know
lower and upper bounds on $\rk(\lambda')$, and on $j'$.  
\smallskip

\item[(2)]

\noindent {\sf Question:} Is there an affine semigroup $\Lambda$ which is Cohen-Macaulay over
some field $k$ but not Cohen-Macaulay over some other field $k'$ ?
Is there an affine semigroup $\Lambda$ which is Cohen-Macaulay over
some field $k$ but not homotopically Cohen-Macaulay ?

Moreover, in case the answer to the first question is yes, it is 
interesting to know whether either or both of the sets of characteristics
for which $\Lambda$ is or is not Cohen-Macaulay can be infinite. 
\smallskip

\item[(3)]
Work of Conca, Herzog, Trung and Valla \cite{CHTV} shows that
for every $\Lambda$ and field $k$ there exists $r>0$ such that 
the rank-selected subposet
$\Lambda_r =\{\la\in\Lambda\mid r \text{ divides } \rk(\la)\}$
has the property that all lower intervals  $(\nnull,\la)$
in $\Lambda_r$ are Cohen-Macaulay over $k$.

\noindent{\sf Question:}  Does this have a combinatorial explanation?

There is an analogous result for bigraded affine semigroups.
Let $f : \Lambda \rightarrow \NN^2$ be a  map of semigroups such that
$f^{-1}((0,1)) \cup f^{-1}((1,0))$ is a generating set of $\Lambda$ 
consisting only of elements of rank $1$.
For $\gamma' \in \NN^2$ denote by $\Lambda_{\gamma'}$ the affine semigroup
of all $\lambda \in \Lambda$ such that $f(\lambda)$ is a multiple of $\gamma'$.
Then (by \cite{CHTV}) there is a $\gamma \in \NN^2$ such that
for all $\gamma' \geq \gamma$ -- this order relation is taken in $\NN^2$ -- 
$\Lambda_{\gamma'}$ is Cohen-Macaulay over $k$.
Of course, in general $k[\Lambda_{\gamma'}]$ does not even have to
be homogeneous. 

\noindent{\sf Question:}  Does this have a combinatorial explanation?
Is there a version of this result for the property ``homotopically 
Cohen-Macaulay''? 
\smallskip

\item[(4)]
Let again $f : \Lambda \rightarrow \NN^2$ be a  map of semigroups such that
$f^{-1}((0,1)) \cup f^{-1}((1,0))$ is a generating set of $\Lambda$ 
consisting only of elements of rank $1$.
A result by Blum \cite{Bl} says that if $k[\Lambda]$ is Koszul then
for $\gamma \in \NN^2$ the affine semigroup ring $k[\Lambda_\gamma]$ is
Koszul as well. (Here we use the notation from Problem (3).) 

\noindent {\sf Question:} Does this have a combinatorial explanation?

Let us see how this result relates to weighted Segre products. 
For weighted Segre products we are given maps $\rk : \Lambda_1 \rightarrow \NN$
and $g : \Lambda_2 \rightarrow \NN$ such that $\rk$ is the rank function of
$\Lambda_1$ and $g$ is a strictly monotone map of affine semigroups.   
These two together give a map $f := (\rk,g) : \Lambda_1 \times \Lambda_2
\rightarrow \NN^2$. Let $\Gamma = f(\Lambda_1 \times \Lambda_2) \cap \{ (a,a) | a \in \NN\}$.
Then $\Lambda_1 \circ_g \Lambda_2$ is the affine semigroup
$(\Lambda_1 \times \Lambda_2)_\Gamma := f^{-1}(\Gamma)$.
If $\Gamma$ is generated by a single element $\gamma$ then $(\Lambda_1 \times \Lambda_2)_\Gamma$
is a diagonal in the sense of \cite{CHTV} and \cite{Bl}
---  except that $f$ usually does not fulfill that  
$f^{-1}((0,1)) \cup f^{-1}((1,0))$ is a generating set of $\Lambda_1 \times
\Lambda_2$.

We can also interpret this result as a result about the weighted Segre product of
two affine semigroup rings $k[\Lambda]$ and $k[\Gamma]$. Suppose that
$k[\Lambda]$ is standard bigraded (i.e. graded with grading in $\NN^2$ and
generated by elements of degree $(1,0)$, and $(0,1)$) and $k[\Gamma]$ is
bigraded in $\NN^2$ by some grading $g$. Then the obvious extension of
the symbol $k[\Lambda] \circ_g k[\Gamma]$ to this situation gives 
$k[\Lambda_\gamma]$ if we take $k[\Gamma] = k[t]$ and grade $t$ by $\Gamma$.

\noindent {\sf Question:} Is there a result about the 
preservation of Koszulness for
weighted Segre products of bigraded Koszul algebras ?

\noindent {\sf Question:} Is there a result about the preservation of Cohen-Macaulayness
for weighted Segre products of bigraded posets ?
\smallskip

\item[(5)]
Let $\Lambda_d$ denote the affine semigroup generated by all vectors
$\la=(\la_1, \ldots, \la_d)\in \NN^d$ such that $\sum \la_i =d$,
except $(1,\ldots,1)$. It has been shown for $d=3$ that 
$\Lambda_d$ is Cohen-Macaulay, or equivalently that $k[\Lambda_d]$ is
Koszul \cite{Cav}. For $d \geq 4$ the question is still open.

\noindent {\sf Question:} Is $\Lambda_d$ Cohen-Macaulay for all $d$ ?
\smallskip

\item[(6)] Define the Rees product $R_n$ as in Example \ref{deranged}.
Being homotopically Cohen-Macaulay we know that $R_n$ is homotopy
equivalent to a wedge of spheres of dimension 
$n-1$. We conjecture that the number of
spheres in this wedge is the derangement number $D_n$,
i.e., the number of permutations 
in the symmetric group $S_n$ without fixed points.
For $n \leq 7$ we have verified by computer that the homology of $R_n$ is
concentrated in top dimension and is free of rank $D_n$.
Since the poset is homotopically Cohen-Macaulay this implies the conjecture
for $n \leq 7$.

The evidence for this conjecture, 
other than computation for small cases, is a natural relationship 
with another poset, which has already been seen to have that
homotopy type. Namely, let $K_n$ be the set of words of pairwise distinct letters over
$[n]$ (i.e. an element of $K_n$ is a sequence $a_1\cdots a_k$ where 
$a_i \in [n]$ and $a_i \neq a_j$ for $1 \leq i < j \leq k$).
We order $K_n$ by subword order: $a_1 \cdots a_k \leq b_1 \ldots b_l$ if and only if
there are indices $1 \leq i_1 < \cdots < i_k \leq l$ such that
$a_1 \cdots a_k = b_{i_1} \cdots b_{i_k}$. 
By results of \cite{Fa} and \cite{BjWa2} it follows that $K_n$ is
homotopy equivalent to a wedge of $D_n$ spheres of dimension $n-1$.
The two posets are related by the 
poset map $\phi : K_n \rightarrow R_n$ which sends
$a_1 \cdots a_k$ to $(\{a_1 , \ldots , a_k\},j)$, where $j-1$ is
the number of descents in $a_1 \cdots a_k$. 
{\em Does this map relate the two posets homotopically?}

For $(A,i) \in R_n$ the lower fiber $\phi^{-1} ((R_n)_{\leq (A,i)})$
is the the order ideal $I_{A,i}$ generated by all words which use all 
letters in $A$ and have $i-1$ descents.
This ideal is, as examples 
show (see Table 1), in general {\em not} contractible. However, it
seems to have reduced Euler-characteristic $0$ -- which would suffice since
both our posets are Cohen-Macaulay. 
Clearly, $I_{A,i}$ as a poset only depends on $i$ and the cardinality of $A$. Thus it
suffices to consider the case $A = [n]$. 
The following table lists the homology groups $\widetilde{H}_*(I_{[n],i},{\bf Z})$.
Note that we only list explicitly those homology groups that are non-zero.

\medskip

\begin{tabular}{||l||c|c|c|c|c|c||}
\hline
\hline
$n \backslash i$ & 1 & 2 & 3 & 4 & 5 & 6  \\
\hline
\hline
      &   &   &   &   &   &   \\
1     & 0 &   &   &   &   &   \\
      &   &   &   &   &   &   \\
\hline
      &   &   &   &   &   &    \\
2     & 0 & 0 &   &   &   &    \\
      &   &   &   &   &   &    \\
\hline
      &   &   &   &   &   &    \\
3     & 0 & $\widetilde{H}_1 = \widetilde{H}_2 = {\bf Z}$ & 0 &   &   &    \\
      &   &   &   &   &   &    \\
\hline
      &   &   &   &   &   &    \\
4     & 0 & 0 & 0 & 0 &   &    \\
      &   &   &   &   &   &    \\
\hline
      &   &   &   &   &   &    \\
5     & 0 & $\widetilde{H}_3 = \widetilde{H}_4 = {\bf Z}$ & $\widetilde{H}_3 = \widetilde{H}_4 = {\bf Z}^6$ & $\widetilde{H}_3 = \widetilde{H}_4 = {\bf Z}$ &  0 &    \\
      &   &   &   &   &   &    \\
\hline
      &   &   &   &   &   &    \\
6     & 0 & 0 & $\widetilde{H}_4 = \widetilde{H}_5 = {\bf Z}^{13}$ & $\widetilde{H}_4 = \widetilde{H}_5 = {\bf Z}^{13}$ & 0 & 0 \\
      &   &   &   &   &   &   \\
\hline 
\hline
\end{tabular}

\medskip

\centerline{{\bf Table 1.} Homology groups of $I_{[n],i}$.}
\end{enumerate}
\bigskip

\section{Acknowledgment}

The authors would like to thank J\"urgen Herzog for helpful hints and discussions.

\bibliographystyle{amsplain}

\begin{thebibliography}{xxxx}
\bibitem{AP} Avramov, L.; Peeva, I.: Finite regularity and Koszul algebras.
Amer. J. Math. {\bf 123} (2002), 275--281.
\bibitem{BF} Backelin, J.; Fr\"oberg, R.: Koszul algebras, Veronese subrings
and rings with linear resolutions. Rev. Roumaine Math. Pures Appl.  {\bf 30} (1985), no. 2,
85--97.
\bibitem{B} Baclawski, K.: Cohen-Macaulay ordered sets. J. Algebra {\bf 63} (1980), 226--258.
\bibitem{Bj} Bj\"orner, A.: Topological Methods. In: Handbook of Combinatorics,
eds. R. Graham, M. Gr\"otschel, L. Lov\'asz. North-Holland, Amsterdam,
1995, pp. 1819--1872.
Advances Math. {\bf 43} (1982), 87--100.
\bibitem{BjWa2} Bj\"orner, A.; Wachs, M.: On lexicographically shellable posets.
Trans. Amer. Math. Soc. {\bf 277 } (1983), 323--341.
\bibitem{BWW} Bj\"orner, A.; Wachs, M.; Welker, V.: Poset fiber theorems. 
Trans. Amer. Math. Soc. {\bf xx} (2004), --. (Preprint 2002)
\bibitem{Bl} Blum, S.: Subalgebras of bigraded Koszul Algebras. J. Algebra {\bf 242} (2001),795--809.
\bibitem{CSV} Carlitz, L.; Scoville, R.; Vaughan, T.: Enumeration of pairs of
              permutations and sequences. Bull. Amer. Math. Soc. {\bf 80} (1974), 881--884.
\bibitem{Cav} Caviglia, G.: The pinched Veronese is Koszul, (Preprint 2003).
\bibitem{CHTV} Conca, A.; Herzog, J.; Trung, N. V.; Valla, G.:
Diagonal subalgebras of bigraded algebras and embeddings of blow-ups of projective spaces.
Amer. J. Math. {\bf 119} (1997), no. 4, 859--901.
\bibitem{C} Crona, K.: A new class of Koszul algebras. C. R. Acad. Sci. Paris, Ser. I 
{\bf 323} (1996), 705--710.
\bibitem{Fa} Farmer, F.D.: Cellular homology of posets. Math. Japon. {\bf 23} (1978/79), 607--613.
\bibitem{F} Fr\"oberg, R.: Koszul algebras. In ``Advances in commutative 
ring theory'' (Fez, 1997), 337--350, 
Lecture Notes in Pure and Appl. Math., 205, Dekker, 1999. 
\bibitem{FH} Fr\"oberg, R.; Hoa, L.T.: Segre products and Rees algebras of
face rings. Comm. Alg. {\bf 20} (1992), 3369--3380.
\bibitem{LS} Laudal, A. O.; Sletsj\o e, A.: Betti numbers of monoid algebras.
Applications to $2$-dimensional torus embeddings.
Math. Scand.  {\bf 56} (1985), 145--162.
\bibitem{PRS} Peeva, I.; Reiner, V.; Sturmfels, B.: How to shell a
monoid. Math. Ann. {\bf 310} (1998), no. 2, 379--393. 
\bibitem{Q} Quillen, D.: Homotopy properties of the poset of non-trivial
$p$-subgroups of a finite group. Advances Math. {\bf 28} (1978), 101--128.
\bibitem{St2} Stanley, R.: Binomial posets, M\"obius inversion, and permutation enumeration.
J. Combin. Theory, Ser. A {\bf 20} (1976), 336--356.
\bibitem{St1} Stanley, R.: Enumerative Combinatorics, Vol. 1.
Wadsworth, Monterey, 1986. Second printing, Cambridge Univ. Press, 1997.
\bibitem{St} Stanley, R.: Combinatorics and Commutative Algebra.
Second Edition. Birkh\"auser, Boston, 1995.
\bibitem{Wa} Walker, J.W.: Canonical homeomorphisms of posets. 
Europ. J. Combinatorics {\bf 9} (1988), 97--107.
\end{thebibliography}

\bigskip

\end{document}